\documentclass[12pt]{scarticle}

\usepackage{
amsmath,
amsfonts,
amssymb,
amscd,
amsthm, 
calc, 
enumerate, 
rotating,
}

\usepackage[all]{xy}

\usepackage{caption}[2006/03/21]
\DeclareCaptionFormat{nocaps}{#1 #3\par}

\newtheorem{theorem}{Theorem}[section]
\newtheorem{proposition}[theorem]{Proposition}
\newtheorem{lemma}[theorem]{Lemma}
\newtheorem{corollary}[theorem]{Corollary}

\theoremstyle{definition}
\newtheorem{definition}[theorem]{\bf Definition}

\newenvironment{acknowledgments}{\noindent\bf Acknowledgments.\rm}{\nolinebreak\par}

\newtheorem{texample}[theorem]{Example}
\newenvironment{example}{\begin{texample}\rm}
{\nolinebreak\mbox{}\nolinebreak\hfill\raisebox{-.262mm}{$\square$}\end{texample}}

\newenvironment{classification}{\noindent\sc Mathematics Subject
Classification:\rm}{\mbox{}}
\newenvironment{keywords}{\noindent\sc Keywords:\rm}{\mbox{}}

\renewcommand{\leq}{\leqslant}
\renewcommand{\geq}{\geqslant}

\newcommand{\concat}{\hspace*{.2mm}\widehat{\mbox{ }}\hspace*{.2mm}}

\newcommand{\bigset}[1]{\big\{ #1 \big\}}
\newcommand{\T}{\makebox[4mm][c]{$|_T$}}    
\newcommand{\nT}{\makebox[4.9mm][c]{$\not|_T$}}    
\newcommand{\wh}[1]{\widehat{#1}}
 
\newcommand{\LL}{{L}}
\newcommand{\M}{{\mathfrak M}}  
\newcommand{\A}{\mathcal{A}}  
\newcommand{\B}{\mathcal{B}}  
   
\newcommand{\C}{\mathcal{C}} 
\newcommand{\D}{\mathcal{D}} 
\newcommand{\E}{\mathcal{E}} 
\newcommand{\F}{\mathcal{F}}  
\newcommand{\G}{\mathcal{G}}  
\newcommand{\HH}{\mathcal{H}}  
\renewcommand{\P}{\mathcal{P}}  
\newcommand{\X}{\mathcal{X}}  
\newcommand{\Y}{\mathcal{Y}}  
\newcommand{\bfA}{{\bf A}}
\newcommand{\bfB}{{\bf B}}
\newcommand{\bfS}{{\bf S}}
\newcommand{\0}{{\bf 0}}

\newcommand{\join}{+}
\newcommand{\meet}{\times}

\newcommand{\Meet}{{\textstyle \prod}}

\begin{document}

\title{The finite intervals of the Muchnik lattice}

\author{Sebastiaan A. Terwijn\thanks{Institute
of Discrete Mathematics and Geometry,
Technical University of Vienna,
Wiedner Hauptstrasse 8--10/E104,
A-1040 Vienna, Austria,
terwijn@logic.at.
Supported by the Austrian Science Fund FWF under grant P18713-N12.}
}

\date{\normalsize \today}

\maketitle

\begin{abstract} \noindent
We characterize the finite intervals of the Muchnik lattice by proving that 
they are a certain proper subclass of the finite distributive lattices.  
\\[5pt]

\end{abstract}

\section{Introduction} \label{intro}

The Medvedev lattice and the Muchnik lattice are structures from 
computability theory that were originally defined for their 
connections with constructive logic, but that are of 
independent interest as well. Both can be seen as 
generalizations of the Turing degrees,  and for example 
when Muchnik presented his solution to Post's problem he 
phrased it as a result of the Medvedev lattice.  
In Terwijn~\cite{Terwijnta} the structure of the Medvedev 
lattice $\M$ was investigated, and it was proven there that 
the finite intervals of $\M$ are precisely the finite Boolean algebras, 
and that the infinite intervals of $\M$ all have cardinality $2^{2^{\aleph_0}}$ 
(cf.\ Theorem~\ref{extendeddichotomy}). 
It was noted there that this strong dichotomy does not hold 
for the Muchnik lattice $\M_w$, and that there are many more 
possibilities for intervals in $\M_w$, both for the finite  
and for the infinite ones. 
In this paper we characterize the finite intervals of $\M_w$ 
by proving that they are a certain subclass of the 
finite distributive lattices that can be described using 
elementary lattice theory. 
In the rest of this section we will repeat the necessary 
definitions and some further preliminaries.

The Medvedev lattice, introduced by Medvedev~\cite{Medvedev}, 
is a particular way of specifying Kolmogorov idea of a calculus of problems. 
Let $\omega$ denote the natural numbers and let $\omega^\omega$ be the 
set of all functions from $\omega$ to $\omega$ (Baire space). 
A {\em mass problem\/} is a subset of $\omega^\omega$. 
Every mass problem is associated with the ``problem'' of 
producing an element of it. 
A mass problem $\A$ {\em Medvedev reduces\/} to 
mass problem $\B$, denoted $\A\leq_M \B$,  
if there is a partial computable functional  
$\Psi:\omega^\omega\rightarrow\omega^\omega$ defined on all of $\B$
such that $\Psi(\B)\subseteq \A$.
That is, $\Psi$ is a uniformly effective method for     
transforming solutions to $\B$ into solutions to~$\A$.   
The relation $\leq_M$ induces an equivalence relation on   
mass problems: $\A\equiv_M\B$ if $\A\leq_M\B$ and $\B\leq_M\A$. 
The equivalence class of $\A$ is denoted by $\deg_M(\A)$ and is 
called the {\em Medvedev degree\/} of~$\A$. 
We denote Medvedev degrees by boldface symbols. 
There is a smallest Medvedev degree, denoted by~$\0$, 
namely the degree of any mass problem containing a 
computable function, and there is a largest degree~$\bf 1$, 
the degree of the empty mass problem, of which it is 
absolutely impossible to produce an element. 
A meet operator $\meet$ and a join operator $\join$ are defined 
on mass problems as follows:
For functions $f$ and $g$, as usual define the function $f\oplus g$
by $f\oplus g(2x)=f(x)$ and $f\oplus g(2x+1)= g(x)$. 
Let $n\concat \A = \{ n\concat f: f\in \A \}$, where $\concat$ 
denotes concatenation. 
Define 
$$
\A \join\B = \big \{ f\oplus g: f\in\A \wedge g\in B\big \}
$$
and 
$$
\A \meet\B = 0 \,\widehat{\mbox{ }}\A \cup 1\concat\B. 
$$
The structure $\M$ of all Medvedev degrees, ordered by $\leq_M$ 
and together with $\join$ and $\meet$ is a distributive lattice
(Medvedev~\cite{Medvedev}).

The Muchnik lattice, introduced by Muchnik~\cite{Muchnik}, 
is a nonuniform variant of the Medvedev lattice. 
It is the structure $\M_w$ resulting from the reduction relation on 
mass problems defined by 
$$
\A \leq_w \B \equiv  (\forall f\in \B)(\exists g\in \A)[g\leq_T f]. 
$$
(The ``w'' stands for ``weak''.)
That is, every solution to the mass problem $\B$ can compute a solution 
to the mass problem $\A$, but maybe not in a uniform way. 
It is easy to check that $\M_w$ is a distributive lattice in the same 
way that $\M$ is, with the same lattice operations and $\bf 0$ and~$\bf 1$. 
Notice that $\A\meet \B$ in $\M_w$ simplifies to $\A\cup \B$.

An M-degree is a {\em Muchnik degree\/} if it contains 
a mass problem that is upwards closed under Turing reducibility 
$\leq_T$. The Muchnik degrees of $\M$ form a substructure that is 
isomorphic to~$\M_w$. 
For any mass problem $\A$, let $C(\A)$ denote the upward 
closure of~$\A$ under~$\leq_T$.  
We have the following embeddings:
$$
\begin{array}{ccccc}
\D_T & \hookrightarrow & \M_w & \hookrightarrow & \M \\
\mbox{Turing degrees} &  & \mbox{Muchnik degrees} & & \mbox{Medvedev degrees} \\
& & & & \\
\deg_T(f) & \longmapsto  & \deg_w(\{f\}) & & \\
& & & & \\
&  & \deg_w(\A) & \longmapsto  & \deg_M(C(\A))
\end{array}
$$
More discussion about the elementary properties of $\M$ and $\M_w$ 
can be found in in Rogers' textbook~\cite{Rogers} and       
the survey paper by Sorbi~\cite{Sorbi}.   
Previous results about embeddings of lattices and algebras 
into $\M$ and $\M_w$ can be found in 
Sorbi~\cite{Sorbi1990,Sorbi1991a}.
Binns and Simpson~\cite{BinnsSimpson} contains results 
about lattice embeddings into the lattice of 
$\Pi^0_1$-classes under $\leq_M$ and~$\leq_w$.

Our notation is mostly standard and follows Odifreddi~\cite{Odifreddi}.  
$\Phi_e$ is the $e$-th partial computable functional. 
For countable sets $I\subseteq \omega$ and mass problems 
$\A_i$, $i\in I$, we have the meet operator 
$$
\Meet_{i\in I} \A_i = \bigset{i\,\concat f: i\in I \wedge f\in\A_i}. 
$$
Note that for finite $I$ this is M-equivalent to an iteration 
of the meet operator $\meet$. 
If $a\leq b$ in some partial order, we use the interval notation 
$[a,b] = \{x: a\leq x\leq b\}$. Similarly $(a,b)$ denotes an interval 
without endpoints, and $(a]$ denotes the set $\{x: x\leq a\}$. 
We say that $b$ {\em covers\/} $a$ if $b>a$ and 
there is no $x$ with $a<x<b$.

In the final section of \cite{Terwijnta} some consequences 
of the results of that paper for the Muchnik lattice $\M_w$ were listed. 
Some of these consequences were: 
\begin{itemize}

\item In contrast to $\M$, the lattice $\M_w$ contains nonempty linear 
intervals. 

\item Every finite Boolean algebra is isomorphic to an interval of $\M_w$.  

\item 
Every Muchnik degree is as large as 
set-theoretically possible:
For every mass problem $\F\subseteq 2^\omega$ we have 
$|\deg_{M_w}(\F)|= 2^{2^{\aleph_0}}$. 

\item 
Whereas in $\M$ only countable Boolean algebras can be embedded, 
the dual of $\P(2^\omega)$ is embeddable into $\M_w$ as a Boolean algebra.    

\end{itemize}

A Medvedev degree is a {\em degree of solvability\/} if it contains a
singleton mass problem. 
A mass problem is called {\em nonsolvable\/} if its M-degree
is not a degree of solvability. 
For every degree of solvability $\bfS$ there is a unique minimal M-degree 
$>\bfS$ that is denoted by $\bfS'$ (cf.\ Medvedev \cite{Medvedev}). 
If $\bfS = \deg_M(\{f\})$ then $\bfS'$ is the degree of the mass problem  
$$
\{f\}' = \bigset{n\concat g: f <_T g \wedge \Phi_n(g) = f }. 
$$
Note that $\{\emptyset\}'$ is M-equivalent to the set of all 
noncomputable functions. We will also denote this set by $0'$. 
Note further that $\{f\}' \equiv_w \{g\in \omega^\omega : f<_T g \}$ 
so that in $\M_w$ we can use this simplified version of $\{f\}'$. 
Dyment~\cite{Dyment} proved that 
the degrees of solvability are precisely characterized
by the existence of such an~$\bfS'$. Namely, 
the degrees of solvability are first-order definable 
(both in $\M$ and in $\M_w$) by the formula 
$$
\phi(x) = \exists y \big(x<y \wedge \forall z(x<z \rightarrow y\leq z)\big).
$$
Thus the Turing degrees form a first-order definable substructure 
of both $\M$ and~$\M_w$.  
This has many immediate corollaries, for example that 
the first-order theories of the structures 
$(\M, \leq_M)$ and $(\M_w, \leq_w)$ are undecidable.

\begin{theorem} {\rm (Dyment \cite{Dyment}, 
cf.\ \cite[Theorem 2.5]{Terwijnta})} \label{Dyment}
For Medvedev degrees $\bfA$ and $\bfB$ with $\bfA <_M \bfB$ it holds that
$(\bfA,\bfB)= \emptyset$ if and only if there is a degree of solvability $\bfS$
such that $\bfA \equiv_M \bfB \meet \bfS$, $\bfB\not\leq_M \bfS$, and $\bfB\leq_M \bfS'$.
\end{theorem}

\noindent 
Theorem~\ref{Dyment} also holds for $\M_w$, with a much easier proof. 
We will include a proof here, as a warm-up for Section~\ref{finite}. 

\begin{lemma} \label{subst}
Suppose that $\A$ and $\B$ satisfy 
\begin{equation}   \label{new}
\forall \C \subseteq \A \mbox{ finite } \; ( \B\meet \C \not\leq_w \A ). 
\end{equation}
Then there exists $\C \geq_w \A$ such that $\C \not\geq_w \B$ and 
$\B\meet \C \not\leq_w \A$.  
If moreover $\A \leq_w \B$ then the interval $(\A, \B)$ is infinite. 
\end{lemma}
\begin{proof} 
Since from (\ref{new}) it follows that $\B\not\leq_w \A$, there is 
$f\in \A$ such that $\{f\} \not\geq_w \B$. 
Again by (\ref{new}) we have that $\B\meet \{f\} \not\leq_w \A$, 
so we can take $\C = \{f\}$. 
 
If in addition $\A \leq_w \B$ then we have 
$\A <_w \B \meet \{f\} <_w \B$. 
Since $\A$ and $\B \meet \{f\}$ also satisfy (\ref{new}) we can 
by iteration of the first part of the lemma obtain an infinite 
downward chain in $(\A, \B)$. 
\end{proof}

\begin{theorem} {\rm (Dyments Theorem for $\M_w$)} \label{DymentMw}
For Muchnik degrees $\bfA$ and $\bfB$ with $\bfA <_w \bfB$ it holds that
$(\bfA,\bfB)=\emptyset$ if and only if there is a degree of solvability $\bfS$
such that $\bfA \equiv_w \bfB \meet \bfS$, $\bfB\not\leq_w \bfS$, and $\bfB\leq_w \bfS'$.
\end{theorem}
\begin{proof}
(If) Suppose that $\bfS = \deg_w(\{f\})$ is as in the theorem and  
suppose that $\A\in \bfA$, $\B\in \bfB$, and 
$\B\meet \{f\}\leq_M \C \leq_M \B$. 
If $\C$ does not contain any element of Turing degree $\deg_T(f)$ 
then it follows that $\C \geq_w \B \meet \{f\}'$, because the elements 
of $\C$ that get sent to the $\{f\}$-side are all strictly 
above $f$, hence included in $\{f\}'$. So in this case $\C\geq_w \B$ 
by $\{f\}' \geq_w \B$.

Otherwise $\C$ contains an element of Turing degree $\deg_T(f)$, 
and consequently $\C\leq_w \{f\}$. 
Hence $\C \leq_w \B\meet \{f\} \equiv_w  \A$.

(Only if) Suppose that $(\A,\B)=\emptyset$. 
Then by Lemma~\ref{subst},  
$\A$ and $\B$ do not satisfy condition~(\ref{new}),   
hence there is a finite set $\C\subseteq \A$ such that 
$\B\meet \C\leq_w \A$. There is also an $f\in\C$ such that 
$\{f\}\not\geq_w \B$, for otherwise we would have $\A\geq_w \B$.
Because the interval is empty and $\A \leq_w \B\meet \{f\} <_w \B$
we must have  $\A\equiv_w \B \meet \{f\}$ 
since there is no other possibility for $\B\meet \{f\}$. 
We also have $\B\meet \{f\}' \not\leq_w \A$ because 
both $\{f\}\not\geq_w \B$ and $\{f\}\not\geq_w \{f\}'$. 
Hence $\B\meet \{f\}' \equiv_w \B$, again by emptiness of the 
interval, and in particular $\{f\}' \geq_w \B$. 
So we can take $\bfS$ to be $\deg_w(\{f\})$. 
\end{proof}

\noindent
Let $f$ and $g$ be T-incomparable. Then the interval 
$\big[\{f,g\},\{f\}'\meet \{g\}'\big]$ contains exactly two 
intermediate elements, cf.\ Figure~\ref{fig}. 
\begin{figure}[htb] 
\begin{center} 
\setlength{\unitlength}{1mm}
\begin{picture}(24, 29)
\put(12,0){\makebox[0cm][c]{$\{f, g\}$}} 
\put(-6,12){\makebox[-.5cm][c]{$\{f\}\meet \{g\}'$}}
\put(-6,16){\line(2,1){12}}   
\put(-6,10){\line(2,-1){12}}
\put(30,12){\makebox[.5cm][c]{$\{g\}\meet \{f\}'$}}
\put(30,16){\line(-2,1){12}}
\put(30,10){\line(-2,-1){12}}
\put(12,24){\makebox[0cm][c]{$\{f\}'\meet \{g\}'$}}
\end{picture}
\end{center}
\captionsetup{format=nocaps}
\caption{\label{fig}}
\end{figure} 
This can be generalized to obtain finite intervals of size $2^n$ for any $n$ 
as follows:

\begin{theorem} \label{exactlyn} 
Let $\B$ be any mass problem. 
Let $n\geq 1$ and let $f_1,\ldots,f_n \in\omega^\omega$ be T-incomparable 
such that $\{f_i\}\not\geq_w \B$ for every $i$. 
Then the interval 
$$
\big[\B\meet\{f_1,\ldots,f_n\}, \B\meet\{f_1\}'\meet\ldots\meet\{f_n\}'\big]
$$ 
in $\M_w$ is isomorphic to the finite Boolean algebra $\mathfrak{2}^n$. 
\end{theorem}
\begin{proof}
This was proved in \cite{Terwijnta} for $\M$. It holds for 
$\M_w$ with the same proof. 
\end{proof}

\noindent
Platek~\cite{Platek} proved that $\M$ has the (for a collection  
of sets of reals maximal possible) cardinality $2^{2^{\aleph_0}}$ 
by showing that $\M$ has antichains of that cardinality.  
(He mentions that the result was noted independently by 
Elisabeth Jockusch and John Stillwell.)
In fact, in $\M$ such large antichains occur in every infinite interval:

\begin{theorem}  {\rm (Terwijn~\cite{Terwijnta}) }  \label{extendeddichotomy} 
Let $[\bfA,\bfB]$ be an interval in $\M$ with $\bfA <_M \bfB$. 
Then either 
$[\bfA,\bfB]$ is isomorphic to the Boolean algebra $\mathfrak{2}^n$ 
for some $n\geq 1$, 
or $[\bfA,\bfB]$ contains an antichain of size $2^{2^{\aleph_0}}$. 
In the latter case, assuming CH, it also contains 
a chain of size $2^{2^{\aleph_0}}$. 
\end{theorem}

\noindent
In particular, $\M$'s version of Theorem~\ref{exactlyn} is the {\em only\/} 
way to generate finite intervals of~$\M$. 
As we will see in what follows, the situation for $\M_w$ is 
rather different.

\section{More on chains and antichains}  

Although every countable linear order can be embedded into $\M_w$ 
(because by Lachlan (cf.\ \cite[p529]{Odifreddi}) this already holds 
for the Turing degrees), 
the following result shows that not every countable linear 
order is {\em isomorphic\/} to an interval in~$\M_w$. 
(From Theorem~\ref{main} it will follow that every finite linear order 
is isomorphic to an interval in~$\M_w$.)

\begin{proposition} 
Not every countable linear order is isomorphic to an interval in $\M_w$. 
\end{proposition}
\begin{proof}
Consider the linear order $\omega + \omega^*$ (that is, a copy of $\omega$ 
followed by a reverse copy of~$\omega$).  
Suppose that $\A_n$ and $\B_n$, $n\in \omega$, are mass problems 
such that for all $n$ and $m$
$$
\A_n <_w \A_{n+1} <_w \B_{m+1} <_w \B_m. 
$$
Let $\C = \Meet_{m \in \omega} \B_m \equiv_w \bigcup_{m \in \omega} \B_m$. 
Then for all $n$, $\A_n <_w \C <_w \B_n$, 
so the interval $[\A_0,\B_0]$ is not isomorphic to $\omega + \omega^*$. 
\end{proof}

\begin{proposition} \label{downwardchain}
$\M_w$ contains linear intervals that are countably infinite.
\end{proposition}
\begin{proof}
By Lachlan, cf.\ \cite[p529]{Odifreddi}, every countable distributive lattice
with a least element is embeddable as an initial segment in the Turing degrees.
Consider the linear order $1 + \omega^*$ (a least element plus a reverse
copy of $\omega$) and embed this in the T-degrees:
Let $\F = \{f_n: n\in\omega\}$ be such that $f_{n+1}<_T f_n$ and such that    
$$
h\leq_T f_0 \rightarrow h \mbox{ computable } \vee \exists n \; h \equiv_T f_n
$$  
for every~$h$. 
Let $\B = \{h : h \not\leq_T f_0\}$ and $\A = \B \meet \F$. 
(Note that $\A$ is in fact w-equivalent to~$0'$.) 
Now if $\C \in [\A,\B]$ then $\C$ can be split in disjoint parts
$\C_0$ and $\C_1$ such that $\C_0$ is maximal with the property $\B \leq_w \C_0$
and $\{f: f\leq_T f_0\} \leq_w \C_1$.
Then $\C_1 \subseteq \F$ and $\C \equiv_w \B \meet \C_1$.
So it suffices to analyze all subclasses of~$\F$:
For every $I\subseteq \omega$ consider
$\C_I = \B\meet\{f_n: n\in I\}$. 
If $I$ is infinite then $\C_I \leq_w \F$, hence $\C_I \equiv_w \A$. 
For $I$ and $J$ finite we have
$\C_I\leq_w \C_J$ whenever $\min I \leq \min J$.
So we see that the interval $(\A,\B)$ contains only the countably
many elements $\B \meet \{f_n\}$, $n\in\omega$.
\end{proof}

\noindent
By Proposition~\ref{downwardchain} there are linear nonempty intervals in $\M_w$.  
This contrasts the situation for $\M$, where by Theorem~\ref{extendeddichotomy} 
all the linear intervals are empty. 
So here we already see that Theorem~\ref{exactlyn} is not the only way anymore 
to generate finite intervals.

$\M_w$ contains antichains of size $2^{2^{\aleph_0}}$, using the same 
argument that Platek used for $\M$,  
(starting with an antichain of size $2^{\aleph_0}$ in the Turing degrees, 
form $2^{2^{\aleph_0}}$ incomparable combinations) 
but Proposition~\ref{downwardchain} shows that 
they do not occur in every infinite interval, as we had for $\M$ 
(cf.\ Theorem~\ref{extendeddichotomy}). 
In fact there are intervals with maximal antichains of every possible size: 

\begin{theorem}  \label{antichainsinMw}  
Each of the following possibilities is realized by some 
interval $[\bfA,\bfB]$ in~$\M_w$:
\begin{enumerate}
\item 
$[\bfA,\bfB]$ contains an antichain of size $n$, but not of size $n+1$, 
\item 
$[\bfA,\bfB]$ contains an antichain of size $\aleph_0$, but no uncountable antichain, 
\item 
$[\bfA,\bfB]$ contains an antichain of size $2^{\aleph_0}$, 
but not of size $2^{2^{\aleph_0}}$, 
\item {\rm (Platek~\cite{Platek}) }
$[\bfA,\bfB]$ contains an antichain of size $2^{2^{\aleph_0}}$. 
\end{enumerate}
\end{theorem} 
\begin{proof}
{\em Ad 1.} This follows from Theorem~\ref{exactlyn}. 

{\em Ad 2.} 
Let $x_0 < x_1 < x_2 < \ldots$ be an increasing chain of elements 
in some lattice and let 
$y_0 > y_1 > y_2 < \ldots$ be a decreasing chain of elements in 
the same lattice such that $x_n \, | \, y_n$ for all $n$. 
Let $\LL$ be the free distributive lattice on these sets of elements with 
an additional least element. 
Then $\LL$ is a countable bottomed distributive lattice, so by 
Lachlan \cite[p529]{Odifreddi}, $\LL$ is embeddable into
the Turing degrees as an initial segment. 
Let $\{f_n : n\in \omega\}$ and $\{g_n : n\in \omega\}$  be 
representatives from the image ${\rm Im}(\LL)$ of $\LL$ corresponding to the 
sequences $x_n$ and~$y_n$ respectively, 
so that $f_i \T g_j$ for all $i$ and~$j$ and such that 
for all~$n$, $f_n <_T f_{n+1}$ and $g_{n+1} <_T g_n$. 
Let 
{\setlength\arraycolsep{2pt}
\begin{eqnarray*} 
\B &=& \{h: \forall n \;\; h \not\leq_T f_n, g_n\} \mbox{ and } \\
\A &=& \B \meet \{f_n, g_n : n\in \omega\}. 
\end{eqnarray*}
}%
Then every $\C\in[\A,\B]$ can be split as 
$\C \equiv_w \C_0 \meet \C_1$,
with $\C_0\subseteq \C$ maximal with the property that $\B \leq_w \C_0$
and $\C_1 \subseteq {\rm Im}(\LL)$. 
Claim: the only elements of ${\rm Im}(\LL)$ that are not in $\B$ are 
the $f_n$, $g_n$, $n\in\omega$. 
To see the claim, note that the nonzero elements of ${\rm Im}(\LL)$ are free 
combinations of the $f_n$ and $g_n$. 
Clearly $\B$ is closed under joins. 
By freeness of $L$ it also easily follows that 
${\rm Im}(\LL)- \{f_n,g_n : n\in\omega\}$ is closed under meets. 
Hence ${\rm Im}(\LL)\cap \B$ is closed under meets and joins, and from this  
it easily follows by induction on the complexity of the elements that every 
element in ${\rm Im}(\LL) - \{f_n,g_n : n\in\omega\}$ is in~$\B$. 
This proves the claim. As a consequence, we have 
(by maximality of $\C_0$) that 
$\C_1 \subseteq \{f_n,g_n : n\in\omega\}$. 
Now $\deg_w(\C)$ is determined by $\deg_w(\C_1)$: One easily checks 
that if $\C,\D \in [\A,\B]$ are split as above as 
$\C \equiv_w \C_0 \meet \C_1$ and 
$\D \equiv_w \D_0 \meet \D_1$ then $\C_1 \equiv_w \D_1$ implies that 
$\C \equiv_w \D$. 
In its turn, $\deg_w(\C_1)$ is determined by the minimal~$n$ (if any) 
such that $f_n \in \C_1$ and by whether $\C_1$ contains infinitely or 
finitely many $g_m$'s, and in the latter case by 
the maximal~$m$ (if any) such that $g_m \in \C_1$. 
So we see that there are only countably many possibilities for the 
degree of $\C_1$, and hence for the degree of $\C$, 
and hence $[\A,\B]$ is countable. 

Now consider the mass problems $\C_n = \B\meet\{f_n, g_n\}. $
Clearly $\C_n \,|_w \, \C_m$ if $n\neq m$. So $[\A,\B]$ is countable 
and contains an infinite antichain. 

{\em Ad 3.} 
Let $\LL$ be a countably infinite distributive lattice with a least element
and an infinite antichain.
By Lachlan \cite[p529]{Odifreddi}, $\LL$ is embeddable into
the Turing degrees as an initial segment.
Let $f_n$, $n\in \omega$, be a set of representatives of
all the degrees in the image of $\LL$.
Consider the interval $[\A, \B]$, where 
$\B = \{h : \forall n \; h\not\leq_T f_n \}$ and 
$\A = \B\meet \{f_n : n\in\omega\}$.
Then $[\A, \B]$ contains an infinite antichain of elements of the
form $\B \meet \{f\}$ because $\LL$ contains a corresponding
infinite antichain. 
For $I \subseteq \omega$ let
$\C_I = \B\meet \bigset{f_n: n\in I}$.
Then for incomparable sets $I$, $J\subseteq \omega$  it holds that
$\C_I \, |_w \, \C_J$.
So $[\A, \B]$ contains an antichain of size $2^{\aleph_0}$.
Now if $\C \in [\A, \B]$ then $\C \equiv_w \C_0 \meet \C_1$,
with $\C_0\subseteq \C$ maximal with the property that $\B \leq_w \C_0$
and $\C_1 \subseteq \A$. So the Muchnik degree $\deg_w(\C)$ of
every $\C \in [\A, \B]$  is determined by a countable set $\C_1$,
hence there are at most $2^{\aleph_0}$ many elements in $[\A, \B]$.

{\em Ad 4.} We can apply Plateks argument to any interval 
that contains an antichain of size $2^{\aleph_0}$ of singletons:  
Suppose that the interval $[\A,\B]$ contains the elements $\B\meet \{f_\alpha\}$, 
$\alpha < 2^\omega$, such that the $f_\alpha$ form an antichain in the Turing 
degrees. For $I \subseteq 2^\omega$ let 
$\C_I = \B\meet \bigset{f_\alpha: \alpha\in I}$. Clearly $\C_I \in [\A, \B]$.
Now for incomparable sets $I$, $J\subseteq 2^\omega$  it holds that 
$\C_I \, |_w \, \C_J$, so it suffices to note that there is an antichain of 
size $2^{2^{\aleph_0}}$ in $\P(2^\omega)$.  (For some general notes on chains 
and antichains we refer to~\cite{Terwijnta}.)  
\end{proof}

\noindent
From the proof of Theorem~\ref{antichainsinMw} we can also read off 
some consequences for chains in $\M_w$: 
\begin{enumerate} 

\item 
By Theorem~\ref{exactlyn} there are intervals 
containing chains of size $n$ but not of size $n+1$, 

\item 
By the proof of item~2., and also Proposition~\ref{downwardchain}, 
there are countable intervals with an infinite chain, 

\item 
The example of an interval given in the proof of item 3.\ contains also a 
chain of size $2^{\aleph_0}$, but not of size $2^{2^{\aleph_0}}$.  
This is because $\P(\omega)$ has a chain of size $2^{\aleph_0}$ 
so the same holds with $\omega$ 
replaced by $\{f_n : n\in\omega\}$. 
A chain in the interval of item~3.\ cannot be bigger since the interval 
itself was of size~$2^{\aleph_0}$. 

\item 
CH implies that $\M_w$ has a chain of size $2^{2^{\aleph_0}}$, 
cf.\ \cite{Terwijnta}.
The conditions for the existence of chains of size $2^{2^{\aleph_0}}$ 
in $\P(2^\omega)$, in $\M$, and in $\M_w$ are the same. 
The consistency of the existence of chains of this size also 
follows from the fact that the dual of $\P(2^\omega)$ is embeddable 
into~$\M_w$, cf.\ \cite{Terwijnta}. 

\end{enumerate}

\section{The finite intervals of $\M_w$} \label{finite}

\begin{theorem} {\rm (Sorbi~\cite{Sorbi1990, Sorbi})} \label{embedding}
A countable distributive lattice with 0,1 is embeddable
into $\M$ (preserving 0 and 1)
if and only if 0 is meet-irreducible and
1 is join-irreducible.
\end{theorem}

\noindent
Sorbi proved Theorem~\ref{embedding} by embedding the (unique)
countable dense Boolean algebra into $\M$.
Since this algebra is embeddable into $\M$ 
it also embeds into $\M_w$. (This is because $\leq_M$ is stronger than 
$\leq_w$ and because for any mass problems $\A$ and $\B$ the 
sets $\A \join \B$ and $\A \meet \B$ are the same in $\M$ and in~$\M_w$.) 
In particular every finite distributive lattice is embeddable 
into~$\M_w$.
In the following we consider lattices that are {\em isomorphic\/}
to an interval of~$\M_w$.
In Theorem~\ref{extendeddichotomy} we saw that for $\M$ these were precisely 
the finite Boolean algebras.
Of course no nondistributive lattice can be isomorphic to an
interval in $\M$ or $\M_w$ since both structures are distributive 
(Medvedev \cite{Medvedev}).
In this section we characterize the finite intervals of 
$\M_w$ as a certain subclass of the finite distributive lattices
(Theorem~\ref{main}). 
We start with some illustrative examples. 

\begin{example}  \label{exdiam}
That the diamond lattice is isomorphic to an interval in $\M_w$ 
was already shown in Theorem~\ref{exactlyn}. 
For later purposes we 
show that this way of obtaining a diamond is essentially the only way. 
Suppose that $[\A,\D]$ is an interval in $\M_w$ containing precisely 
two intermediate elements $\B$ and $\C$, 
and that $\B$ and $\C$ are incomparable, cf.\ Figure~\ref{fff}.
\begin{figure}[htb]
\mbox{} \hfill
\xymatrix@R=20pt@C=20pt{
  & \D & \\
\B \ar@{-}[ur] & 
& \C \ar@{-}[ul] \\
& \A \ar@{-}[ur] \ar@{-}[ul] & \\ 
}
\hfill \mbox{} 
\captionsetup{format=nocaps}
\caption{\label{fff}}
\end{figure}
Then $\A$ and $\D$ do not satisfy property~(\ref{new}) of Lemma~\ref{subst} 
because the interval is finite. 
So there is a finite set $\X\subseteq \A$ such that $\D \meet \X \leq_w \A$, 
and hence $\D \meet \X \equiv_w \A$.  
Without loss of generality the elements of $\X$ are pairwise T-incomparable 
and $\{f\} \not\geq_w \D$ for every $f\in \X$. 
Since for every $f\in\X$ we have $\A <_w \D\meet \{f\} <_w \D$ we see that 
$\X$ can contain at most two elements.
If $\X$ would contain only one element $f$ then $[\D\meet\{f\}, \D\meet\{f\}']$ 
would be an initial segment of the interval $[\A,\D]$, which is a contradiction. 
So $\X$ contains precisely two elements, $f_0$ and $f_1$ say. 
Then $[\A,\D]$ contains the interval 
$[\D\meet \{f_0,f_1\}, \D\meet \{f_0\}' \meet \{f_1\}']$, 
which by Theorem~\ref{exactlyn} is isomorphic to the diamond lattice. 
So we must have that $\D \equiv_w \D \meet \{f_0\}' \meet \{f_1\}'$. 
\end{example}

\begin{example}  \label{exdouble}
We continue with Example~\ref{exdiam}. 
We show how to obtain the lattice depicted in Figure~\ref{ff2} 
as an interval. 
\begin{figure}[htb]
\mbox{} \hfill
\xymatrix@R=20pt@C=20pt{
  & & \F \ar@{-}[dl] \ar@{-}[dr] &  \\
  & \D & & \E \ar@{-}[dl] \\
\B \ar@{-}[ur] & & \C \ar@{-}[ul] &\\
& \A \ar@{-}[ur] \ar@{-}[ul] &  & \\ 
}
\hfill \mbox{} 
\captionsetup{format=nocaps}
\caption{\label{ff2}}
\end{figure}
Let $f_0$ and $f_1$ be T-incomparable and let $g_1 \T f_0$ 
be minimal over $f_1$, i.e.\ $g_1 >_T f_1$ and 
there is no function of T-degree strictly between $f_1$ and $g_1$. 
(In lattice theoretic terminology $g_1$ is said to {\em cover\/} $f_1$.) 
Let $\X = \{ h \geq_T f_1 : h \T g_1\}$. 
We then have $\X\meet \{f_1\}' \equiv_w \X\meet \{g_1\}$, 
as is easy to check. 
Define  
{\setlength\arraycolsep{2pt}
\begin{eqnarray*} 
\A &=& \X \meet \{f_0,f_1\}  \\
\B &=& \X \meet \{f_0\}' \meet \{f_1\}   \\
\C &=& \X \meet \{f_0\} \meet \{f_1\}' = \X\meet \{f_0, g_1\} \\
\D &=& \X \meet \{f_0\}' \meet \{f_1\}' = \X\meet \{f_0\}' \meet \{g_1\} \\
\E &=& \X\meet \{f_0\} \meet \{g_1\}'    \\
\F &=& \X\meet \{f_0\}' \meet \{g_1\}' 
\end{eqnarray*}
}%
Then by Example~\ref{exdiam} the intervals $[\A,\D]$ and $[\C,\F]$ 
are both isomorphic to a diamond, hence the whole interval $[\A,\F]$ 
is isomorphic to the configuration of Figure~\ref{ff2}. 
As in Example~\ref{exdiam} one can also show that the above method is essentially 
the {\em only\/} way of obtaining this configuration as an interval of~$\M_w$. 
\end{example}

\noindent
Using similar methods as in the previous examples one can show that 
the lattices from Figures~\ref{fig1} and~\ref{fig2} can be obtained as    
intervals of $\M_w$.  
\begin{figure}[htb]
\parbox[t]{0.45\textwidth}{
\begin{center} 
\setlength{\unitlength}{1mm}
\begin{picture}(0,34.14)
\put(0,0){\circle*{1.5}}
\put(0,0){\line(0,1){14.14}}
\put(-10,24.14){\circle*{1.5}}
\put(10,24.14){\circle*{1.5}}
\put(0,14.14){\circle*{1.5}}
\put(0,34.14){\circle*{1.5}}
\put(0,14.14){\line(-1,1){10}}
\put(0,14.14){\line(1,1){10}}
\put(0,34.14){\line(-1,-1){10}}
\put(0,34.14){\line(1,-1){10}}
\end{picture}
\end{center}
\captionsetup{format=nocaps}
\captionof{figure}{\label{fig1}}
}
\hfill
\parbox[t]{0.45\textwidth}{
\begin{center} 
\setlength{\unitlength}{1mm}
\begin{picture}(0,34.14)
\put(0,34.14){\circle*{1.5}}
\put(0,20){\line(0,1){14.14}}
\put(-10,10){\circle*{1.5}}
\put(10,10){\circle*{1.5}}
\put(0,0){\circle*{1.5}}
\put(0,20){\circle*{1.5}}
\put(0,0){\line(-1,1){10}}
\put(0,0){\line(1,1){10}}
\put(0,20){\line(-1,-1){10}}
\put(0,20){\line(1,-1){10}}

\end{picture}
\end{center}
\captionsetup{format=nocaps}
\captionof{figure}{\label{fig2}}
}
\end{figure}
For the first one uses an $f$ branching into two 
incomparable elements $f_0$ and $f_1$, and 
for the second one uses incomparable elements $f_0$ and $f_1$ 
and their join $f_0\oplus f_1$.

Next we show that not every finite distributive lattice is isomorphic 
to an interval in $\M_w$. 

\begin{proposition} \label{counterexample} 
The double diamond lattice from Figure~\ref{figccc} is not 
isomorphic to an interval in $\M_w$. 
\end{proposition}
\begin{figure}[htb]
\mbox{} \hfill
\xymatrix@R=20pt@C=20pt{
  & \G & \\
\E \ar@{-}[ur] & 
& \F \ar@{-}[ul] \\
& \D \ar@{-}[ur] \ar@{-}[ul] & \\ 
\B \ar@{-}[ur] & 
& \C \ar@{-}[ul] \\
& \A \ar@{-}[ur] \ar@{-}[ul] & \\ 
}
\hfill \mbox{} 
\captionsetup{format=nocaps}
\caption{\label{figccc}}
\end{figure}
\begin{proof}
Assume for a contradiction that the interval $[\A,\G]$ 
is isomorphic to the lattice of Figure~\ref{figccc}. 
As in Example~\ref{exdiam} we can argue that there is a finite 
set $\X \subseteq \A$ such that $\A \equiv_w \X \meet \G$. 
Using the same reasoning as before we can argue that 
$\X$ contains precisely two T-incomparable elements $f_0$ and $f_1$
with $\{f_0\}, \{f_1\}\not\geq_w \G$.  
(If $\X$ contained at least three of such elements then 
by Theorem~\ref{exactlyn} the interval 
$[\A,\G]$ would contain a copy of $\mathfrak{2}^3$, but 
the interval contains only $7$ elements so this is impossible.) 
Since by Example~\ref{exdiam} there is only one way of obtaining a 
diamond, there are T-incomparable $g_0$ and $g_1$ with 
$\{g_0\}, \{g_1\} \not\geq_T \G$ such that 
{\setlength\arraycolsep{2pt}
\begin{eqnarray*} 
\A &\equiv_w& \G \meet \{f_0,f_1\} \\
\B &\equiv_w& \G \meet \{f_0\}' \meet \{f_1\} \\
\C &\equiv_w& \G \meet \{f_0\} \meet \{f_1\}' \\
\D &\equiv_w& \G \meet \{f_0\}' \meet \{f_1\}' 
\equiv_w \G \meet \{g_0,g_1\} \\
\E &\equiv_w& \G \meet \{g_0\}' \meet \{g_1\} \\
\F &\equiv_w& \G \meet \{g_0\} \meet \{g_1\}' \\
\G &\equiv_w& \G \meet \{g_0\}' \meet \{g_1\}'  
\end{eqnarray*}
}%
From the two equations for $\D$ it follows that 
$\{g_0,g_1\} >_w \{f_0,f_1\}$. Now there are two cases: 
\begin{itemize} 

\item Both $g_i$'s are T-above both $f_j$'s. But then we have 
$$
\D \equiv_w \G \meet \{f_0\}' \meet \{f_1\}' \leq_w \G \meet \{f_0\oplus f_1\} 
<_w \G \meet \{g_0,g_1\} \equiv_w \D, 
$$
a contradiction. (The second to last inequality is strict since 
$f_0\oplus f_1 <_T g_0, g_1$ because $g_0$ and $g_1$ are incomparable.) 

\item The $g_i$'s are not both above $f_0$ and $f_1$. Hence either 
there is precisely one $g_i$ above each $f_j$, or there are 
precisely two $g_i$'s above one $f_j$. 
In both cases there is at least one $g_i$ T-incomparable to an $f_j$, 
say that $f_0 \T g_1$. 
Now consider the element $\mathcal{H} = \G \meet \{f_0\} \meet \{g_1\}'$. 
Clearly $\mathcal{H} \in [\A,\G]$. But $\mathcal{H} $ is w-incomparable to $\D$: 
$\mathcal{H} \not\leq_w \D$ because $\mathcal{H} \not\leq_w \{g_1\}$, and 
$\D \not\leq_w \mathcal{H} $ because $\D \not\leq_w \{f_0\}$. 
So again we have reached a contradiction, because $[\A,\G]$ does not 
contain an element incomparable to $\D$. 
\end{itemize}
Since both cases are contradictory we conclude that it is impossible 
that $[\A,\G]$ is isomorphic to the double diamond. 
\end{proof}

\noindent
We will see in Theorem~\ref{main} that the double diamond lattice of 
Figure~\ref{figccc} is the smallest possible counterexample.

Let us recall some elementary lattice theory from Gr\"atzer~\cite{Graetzer}. 
Let $L$ be a distributive lattice. 
$J(L)$ denotes the set of all {\em nonzero\/} 
join-irreducible elements of~$L$. 
$J(L)$ is a poset under the ordering of~$L$. 
For $a\in L$ define 
$$
r(a) = \bigset{x \in J(L) : x \leq a }.
$$
For a poset $P$ let $H(P)$ be the collection of downwards closed 
subsets of $P$, partially ordered by inclusion. 
Then $H(P)$ is a distributive lattice, and we have

\begin{theorem} {\rm (\cite[Theorem II.1.9]{Graetzer})} \label{isom}
For any finite distributive lattice $L$ the mapping 
$a \mapsto r(a)$ is an isomorphism between $L$ and $H(J(L))$.
\end{theorem}

\noindent
Thus the mappings $J$ and $H$ are inverses of each other, and they 
relate the class of finite distributive lattices with the class of 
all finite posets.

Say that a lattice $L$ contains another lattice $L'$ as 
a {\em subinterval\/} if there is an interval 
$[a,b] \subseteq L$ such that $[a,b] \cong L'$. 
Note that this is not the same as saying 
that $L'$ is a sublattice of $L$. For example, 
the free distributive lattice on three elements 
$F_D(3)$, depicted in Figure~\ref{FD3}, 
contains the double diamond 
of Figure~\ref{figccc} as a sublattice, but not  
as a subinterval. 
\begin{figure}[htb]
\begin{center}
\setlength{\unitlength}{0.8mm}
\begin{picture}(0, 60)
\put(0,10){\circle*{1.5}}
\put(-15,40){\circle*{1.5}}
\put(15,40){\circle*{1.5}}
\put(0,30){\circle*{1.5}}
\put(0,50){\circle*{1.5}}
\put(0,30){\line(-3,2){15}}
\put(0,30){\line(3,2){15}}
\put(0,50){\line(-3,-2){15}}
\put(0,50){\line(3,-2){15}}
\put(15,20){\circle*{1.5}}
\put(-15,20){\circle*{1.5}}
\put(0,10){\line(-3,2){15}}
\put(0,10){\line(3,2){15}}
\put(0,30){\line(-3,-2){15}}
\put(0,30){\line(3,-2){15}}

\put(0,0){\circle*{1.5}}
\put(15,10){\circle*{1.5}}
\put(-15,10){\circle*{1.5}}
\put(0,20){\circle*{1.5}}
\put(0,0){\line(-3,2){15}}
\put(0,0){\line(3,2){15}}
\put(0,20){\line(-3,-2){15}}
\put(0,20){\line(3,-2){15}}

\put(0,0){\line(0,1){10}}
\put(-15,10){\line(0,1){10}}
\put(15,10){\line(0,1){10}}
\put(0,20){\line(0,1){10}}

\put(0,30){\line(0,1){10}}
\put(-15,40){\line(0,1){10}}
\put(15,40){\line(0,1){10}}
\put(0,50){\line(0,1){10}}

\put(0,40){\circle*{1.5}}
\put(15,50){\circle*{1.5}}
\put(-15,50){\circle*{1.5}}
\put(0,60){\circle*{1.5}}
\put(0,40){\line(-3,2){15}}
\put(0,40){\line(3,2){15}}
\put(0,60){\line(-3,-2){15}}
\put(0,60){\line(3,-2){15}}

\put(-30,30){\circle*{1.5}}
\put(15,30){\circle*{1.5}}
\put(30,30){\circle*{1.5}}
\put(-30,30){\line(3,2){15}}
\put(-30,30){\line(3,-2){15}}
\put(30,30){\line(-3,2){15}}
\put(30,30){\line(-3,-2){15}}
\put(15,30){\line(-3,2){15}}
\put(15,30){\line(-3,-2){15}}

\end{picture}
\captionsetup{format=nocaps}
\caption{\label{FD3}}
\end{center}
\end{figure}

\begin{definition} 
We call a finite distributive lattice {\em double diamond-like\/} if it has 
at least two elements, it has no largest and smallest nonzero join-irreducible 
element, and $0$ is not the meet of two nonzero join-irreducible elements 
one of which is maximal. 
\end{definition}

\begin{figure}[htb]
\begin{center} 
\setlength{\unitlength}{.8mm}
\begin{picture}(110, 70)
\put(-20,0){\circle*{1.5}}
\put(-30,30){\circle*{1.5}}
\put(-10,30){\circle*{1.5}}
\put(-20,20){\circle*{1.5}}
\put(-20,40){\circle*{1.5}}
\put(-20,20){\line(-1,1){10}}
\put(-20,20){\line(1,1){10}}
\put(-20,40){\line(-1,-1){10}}
\put(-20,40){\line(1,-1){10}}
\put(-10,10){\circle*{1.5}}
\put(-30,10){\circle*{1.5}}
\put(-20,0){\line(-1,1){10}}
\put(-20,0){\line(1,1){10}}
\put(-20,20){\line(-1,-1){10}}
\put(-20,20){\line(1,-1){10}}
\put(-20,0){\circle*{1.5}}
\put(-10,10){\circle*{1.5}}
\put(0,40){\circle*{1.5}}
\put(-10,30){\circle*{1.5}}
\put(-10,50){\circle*{1.5}}
\put(-10,30){\line(1,1){10}}
\put(-10,50){\line(-1,-1){10}}  
\put(-10,50){\line(1,-1){10}}
\put(0,20){\circle*{1.5}}
\put(-10,10){\line(1,1){10}}
\put(-10,30){\line(1,-1){10}}
\put(-10,10){\circle*{1.5}}
\put(30,0){\circle*{1.5}}
\put(30,0){\line(0,1){14.14}}
\put(20,24.14){\circle*{1.5}}
\put(40,24.14){\circle*{1.5}}
\put(30,14.14){\circle*{1.5}}
\put(30,34.14){\circle*{1.5}}
\put(30,14.14){\line(-1,1){10}}
\put(30,14.14){\line(1,1){10}}
\put(30,34.14){\line(-1,-1){10}}
\put(30,34.14){\line(1,-1){10}}
\put(50,5){\circle*{1.5}}
\put(50,5){\line(0,1){14.14}}
\put(40,29.14){\circle*{1.5}}
\put(60,29.14){\circle*{1.5}}
\put(50,19.14){\circle*{1.5}}
\put(50,39.14){\circle*{1.5}}
\put(50,19.14){\line(-1,1){10}}
\put(50,19.14){\line(1,1){10}}
\put(50,39.14){\line(-1,-1){10}}
\put(50,39.14){\line(1,-1){10}}
\put(50,53.28){\circle*{1.5}}
\put(50,39.14){\line(0,1){14.14}}
\put(60,34.14){\circle*{1.5}}
\put(80,34.14){\circle*{1.5}}
\put(70,24.14){\circle*{1.5}}
\put(70,44.14){\circle*{1.5}}
\put(70,24.14){\line(-1,1){10}}
\put(70,24.14){\line(1,1){10}}
\put(70,44.14){\line(-1,-1){10}}
\put(70,44.14){\line(1,-1){10}}
\put(70,58.28){\circle*{1.5}}
\put(70,44.14){\line(0,1){14.14}}
\put(30,0){\line(4,1){20}}
\put(20,24.14){\line(4,1){40}}
\put(40,24.14){\line(4,1){40}}
\put(30,14.14){\line(4,1){40}}
\put(30,34.14){\line(4,1){40}}
\put(50,53.28){\line(4,1){20}}
\put(120,0){\circle*{1.5}}
\put(110,30){\circle*{1.5}}
\put(130,30){\circle*{1.5}}
\put(120,20){\circle*{1.5}}
\put(120,40){\circle*{1.5}}
\put(120,20){\line(-1,1){10}}
\put(120,20){\line(1,1){10}}
\put(120,40){\line(-1,-1){10}}
\put(120,40){\line(1,-1){10}}
\put(130,10){\circle*{1.5}}
\put(110,10){\circle*{1.5}}
\put(120,0){\line(-1,1){10}}
\put(120,0){\line(1,1){10}}
\put(120,20){\line(-1,-1){10}}
\put(120,20){\line(1,-1){10}}
\put(120,0){\circle*{1.5}}
\put(130,10){\circle*{1.5}}
\put(140,40){\circle*{1.5}}
\put(130,30){\circle*{1.5}}
\put(130,50){\circle*{1.5}}
\put(130,30){\line(1,1){10}}
\put(130,50){\line(-1,-1){10}}  
\put(130,50){\line(1,-1){10}}
\put(140,20){\circle*{1.5}}
\put(130,10){\line(1,1){10}}
\put(130,30){\line(1,-1){10}}
\put(130,10){\circle*{1.5}}
\put(100,40){\circle*{1.5}}
\put(120,40){\circle*{1.5}}
\put(110,50){\circle*{1.5}}
\put(110,30){\line(-1,1){10}}
\put(110,50){\line(-1,-1){10}}
\put(110,50){\line(1,-1){10}}
\put(120,20){\circle*{1.5}}
\put(100,20){\circle*{1.5}}
\put(110,10){\line(-1,1){10}}
\put(110,30){\line(-1,-1){10}}
\put(110,10){\circle*{1.5}}
\put(120,60){\circle*{1.5}}
\put(120,60){\line(-1,-1){10}}  
\put(120,60){\line(1,-1){10}}
\put(120,20){\circle*{1.5}}
\end{picture}
\end{center}
\caption{Some double diamond-like lattices. \label{more}}
\end{figure}
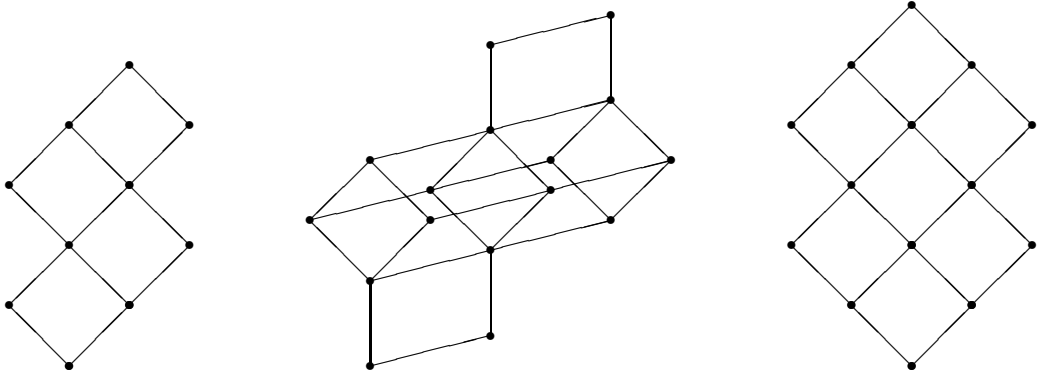

\noindent
The double diamond lattice from Figure~\ref{figccc} is the smallest 
example of a double diamond-like lattice. 
Figure~\ref{more} shows some other examples. 

\begin{theorem} \label{noccc}
For a finite distributive lattice $L$, 
the poset $J(L)$ is an initial segment of an upper semilattice 
if and only if $L$ has no double diamond-like subinterval. 
\end{theorem}
\begin{proof}
Suppose that $J(L)$ is an initial segment of an upper semilattice 
and let $[a,b]$ be an interval in $L$. We prove that 
$[a,b]$ is not double diamond-like. 
To this end, suppose that $a\neq b$ and that $J([a,b])$ has no 
largest and no smallest element. We have to prove that then the~0 of 
$[a,b]$, which is~$a$, is a meet of two elements of $J([a,b])$,  
one of which is maximal. 
By the assumption $a\neq b$ we have $J([a,b])\neq\emptyset$, 
and therefore there are at least two maximal 
nonzero join-irreducible elements $y_0$ and $y_1$ in $[a,b]$ 
and at least two minimal ones, $x_0$ and $x_1$ say. 
Suppose that both $y$'s are above both $x$'s, so that   
$J([a,b])$ contains the configuration of Figure~\ref{JL}. 
\begin{figure}[htb]
\begin{center}
\setlength{\unitlength}{1.5mm}
\begin{picture}(10,13)

\put(0,0){\circle*{1.2}}
\put(0,10){\circle*{1.2}}
\put(10,0){\circle*{1.2}}
\put(10,10){\circle*{1.2}}

\put(0,0){\line(0,1){10}}
\put(10,0){\line(0,1){10}}
\put(0,0){\line(1,1){10}}
\put(0,10){\line(1,-1){10}}

\end{picture}
\captionsetup{format=nocaps}
\caption{\label{JL}}
\end{center}
\end{figure}
We cannot immediately conclude from this that $J(L)$ contains 
the same configuration, for $J([a,b])$ and $J(L)$ can even be disjoint. 
Nevertheless, suppose that $y_0$ is join-reducible in $L$ as 
$y_0 = z_0 \join z_1$, with $z_0 | z_1$. 
By Lemma~\ref{hulp} we can choose $z_0$ and $z_1$ such that 
$z_0\meet x_0 \neq z_0\meet x_1$ and $z_0\not\leq y_1$. 
Then in $L$ the set $\{z_0\meet x_0, z_0\meet x_1, z_0, y_1\}$, 
is partially ordered as in Figure~\ref{JL}. 
Continuing in this way we can reduce the configuration until 
the top element $y_0$ has become join-irreducible, and of course 
we can reduce $y_1$ in the same way. 
Then $L$ contains the configuration of Figure~\ref{JL} with 
both top elements in $J(L)$. 
But the bottom elements always bound nonzero join-reducible elements, 
so we see that $L$ contains Figure~\ref{JL} with all four elements 
in~$J(L)$. 
But this contradicts that $J(L)$ is an initial of an upper semilattice
and hence that the bottom two elements should have a least 
upper bound in~$J(L)$. 
We conclude that in the original configuration $\{x_0,x_1,y_0,y_1\}$ in 
$J([a,b])$ it is not possible that both $y$'s are above both $x$'s,   
so there is at least one $x$ that is not below a $y$, say 
$x_0 \,|\, y_1$. 
But then for the join-irreducible elements $x_0$ and $y_1$ in    
$[a,b]$ we have $a = x_0 \meet y_1$: If $x_0 \meet y_1> a$ then 
$x_0$ would bound some join-irreducible element $>a$, contradicting 
that $x_0$ is minimal in $J([a,b])$.  
This shows that $[a,b]$ is not double diamond-like.   

Conversely, if $L$ is a lattice such that 
$J(L)$ is not an initial segment of any upper semilattice
then it must contain two incomparable elements $x_0$ and $x_1$ without 
a least upper bound. If $x_0$ and $x_1$ would have no or only one upper 
bound in $J(L)$ then $J(L)$ could be consistently extended to 
an upper semilattice, so there must be at least {\em two\/} 
incomparable minimal upper bounds $y_0$ and $y_1$ for both $x_0$ and $x_1$. 
But then the poset $\{x_0,x_1,y_0,y_1\}$ is isomorphic to the 
configuration in Figure~\ref{JL} (with possibly other elements in 
between the $x$'s and $y$'s). Denote this subposet of $J(L)$ by~$P$.    
Now it is not hard to check that $L$ has a double diamond-like subinterval. 
Namely by Theorem~\ref{isom} we have $L \cong H(J(L))$. 
Consider the interval $[X,Y]$ in $H(J(L))$ defined by 
{\setlength\arraycolsep{2pt}
\begin{eqnarray*} 
X &=& \bigset{x\in J(L): x< x_0 \vee x< x_1},  \\
Y &=& \bigset{x\in J(L): x\leq y_0 \vee x\leq y_1}.
\end{eqnarray*}
}%
Then clearly $[X,Y]$ is double diamond-like, 
because $X\neq Y$ and $[X,Y] = H(P)$ and so 
$J([X,Y])$ $=$ $J(H(P)) \cong P$. 
\end{proof}

\begin{lemma}\label{hulp}
In the proof of Theorem~\ref{noccc} above, if $y_0 = z_0 \join z_1$ in 
$L$, $z_0|z_1$, then we can choose such $z_0$ and $z_1$ with 
$z_0\meet x_0 \neq z_0\meet x_1$ and $z_0\not\leq y_1$. 
\end{lemma}
\begin{proof}
Suppose that $y_0 = z_0 \join z_1$ in $L$, with $z_0 | z_1$. 
Note that $z_0$ and $z_1$ cannot be both in $[a,b]$. 
Suppose that 
\begin{equation} \label{star}
\forall v,w \in L \big( v|w \wedge y_0= v\join w \rightarrow 
v,w \notin [a,b] \big).
\end{equation} 
Consider $z_0$ and $a\join z_1$. If $a\join z_1 = y_0$ then this 
contradicts (\ref{star}) (because both $a, z_1< y_0$ they must be 
incomparable in this case). 
If $a\join z_1 < y_0$ then by $(a\join z_1) \join z_0 = y_0$  
we again contradict (\ref{star}). 
Hence (\ref{star}) is false, and if $y_0 = z_0 \join z_1$ with $z_0 | z_1$ in $L$ 
we can always choose $z_0\notin [a,b]$ and $z_1\in [a,b]$. 
In this case $z_0 \join a = y_0$, for if $z_0 \join a < y_0$ then 
by $(a\join z_0) \join z_1 = y_0$ we would have $y_0$ join-reducible in $[a,b]$, 
contradiction. 
Hence for every $c \in [a,y_0]$ it holds that $z_0\join c = y_0$, 
and in particular 
\begin{equation} \label{star2}
z_0\join x_1 = z_0\join x_1 = y_0.
\end{equation}
Now we also have $z_0 \not\leq y_1$ because otherwise 
$y_0 = z_0\join x_0 \leq y_1$, contradiction. 

Finally we prove that $z_0 \meet x_0 \neq z_0\meet x_1$. 
Suppose that $z_0 \meet x_0 = z_0\meet x_1$. 
Because by (\ref{star2})  it holds that 
$z_0\join x_1 = z_0\join x_1$  we have 
$$
\begin{array}{rcll}
x_1 &=& (z_0 \meet x_1) \join x_1 \\
&=& (z_0 \meet x_0) \join x_1 & \\
&=& (z_0 \join x_1) \meet (x_0 \join x_1)  &
\mbox{ (by distributivity) }\\
&=& (z_0 \join x_0) \meet (x_0 \join x_1)  & \\
&=& x_0 \join (z_0 \meet x_1) \\
&\geq& x_0   
\end{array}
$$
From this contradiction we conclude that $z_0 \meet x_0 \neq z_0\meet x_1$. 
\end{proof}

\begin{example} \label{exmain}
Before giving the general result of how to obtain lattices 
as intervals of $\M_w$ we give a specific example to illustrate 
the method. 
Figure~\ref{procedure} depicts the procedure to obtain a given lattice $L$ 
as an interval of~$\M_w$. 
\begin{figure}[htb]
\begin{center} 
\setlength{\unitlength}{1mm}
\begin{picture}(100,120)
\put(-20,116){$L$}
\put(0,80){\circle*{1}}
\put(-10,110){\circle*{1}}
\put(-10,110){\circle{2.5}}
\put(10,110){\circle*{1}}
\put(0,100){\circle*{1}}
\put(20,100){\circle*{1}}
\put(20,100){\circle{2.5}}
\put(0,120){\circle*{1}}
\put(0,100){\line(-1,1){10}}
\put(0,100){\line(1,1){10}}
\put(20,100){\line(-1,1){10}}
\put(20,100){\line(-1,-1){10}}
\put(0,120){\line(-1,-1){10}}
\put(0,120){\line(1,-1){10}}
\put(10,90){\circle*{1}}
\put(10,90){\circle{2.5}}
\put(-10,90){\circle*{1}}
\put(-10,90){\circle{2.5}}
\put(0,80){\line(-1,1){10}}
\put(0,80){\line(1,1){10}}
\put(0,100){\line(-1,-1){10}}
\put(0,100){\line(1,-1){10}}

\put(40,100){$\Longrightarrow$}

\put(100,116){$I(J(L))$}
\put(100,93){\circle*{1}} \put(102,92){$f_1$}
\put(100,108){\circle*{1}} \put(102,107){$g_1$}
\put(80,93){\circle*{1}}  \put(74.5,92){$f_0$}
\put(80,108){\circle*{1}}  \put(74.5,107){$g_0$}
\put(100,93){\line(0,1){15}}
\put(80,93){\line(0,1){15}}
\put(80,108){\line(4,-3){20}}

\put(90,70){%
\begin{rotate}{-90}
\hspace*{-0.6cm}
\raisebox{-.1cm}{
$\Longrightarrow$
}
\end{rotate}}

\put(100,54){$H(I(J(L)))$}
\put(120,55){\makebox[0cm][r]{%
\newcommand{\ttt}[1]{\makebox[45pt][c]{$#1$}}
\xymatrix@R=20pt@C=-13pt{
  & \ttt{\{f_0,f_1,g_0,g_1\}} \ar@{-}[dl] \ar@{-}[dr]  & &  \\
\ttt{\{f_0,f_1,g_0\}} \ar@{-}[dr] & & 
\ttt{\{f_0,f_1,g_1\}} \ar@{-}[dl] \ar@{-}[dr] &  \\
  & \ttt{\{f_0,f_1\}} & & \ttt{\{f_1,g_1\}} \ar@{-}[dl] \\
\ttt{\{f_0\}} \ar@{-}[ur] & & 
\ttt{\{f_1\}} \ar@{-}[ul] &\\
& \ttt{\emptyset} \ar@{-}[ur] \ar@{-}[ul] &  & \\ 
}
}}

\put(40,27){$\stackrel{F}{\Longleftarrow}$}

\put(-20,52){$F(H)$}   
\put(30,52){\makebox[0cm][r]{%
\xymatrix@R=20pt@C=20pt{
  & \mathcal{H}  \ar@{-}[dl] \ar@{-}[dr]  & &  \\
\G \ar@{-}[dr] & & \F \ar@{-}[dl] \ar@{-}[dr] &  \\
  & \D & & \E \ar@{-}[dl] \\
\B \ar@{-}[ur] & & \C \ar@{-}[ul] &\\
& \A\ar@{-}[ur] \ar@{-}[ul] &  & \\ 
}
}}
\end{picture}
\end{center}
\caption{Procedure to obtain an interval in $\M_w$ isomorphic to 
a given~$L$. \label{procedure}}    
\end{figure}
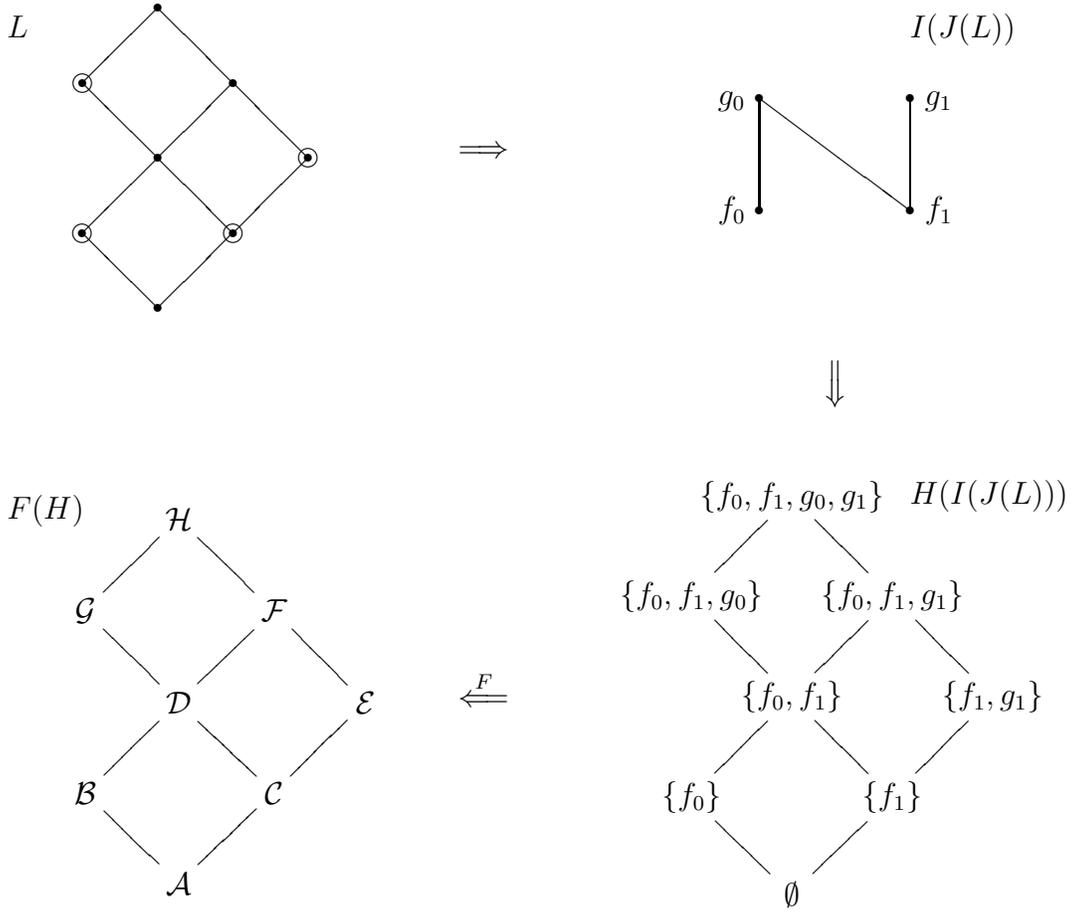
The top left side of the picture shows an example of a finite distributive 
lattice, with its nonzero join-irreducible elements circled.  
The partial order $J(L)$ is depicted on the top right. 
Now for the lattice $L$ in this particular example 
we can map the poset $J(L)$ to an 
isomorphic configuration $I(J(L))$ in~$\D_T$. 
(The picture remains the same, so we drew it only once.) 
This means that the only relations are the ones indicated in the 
picture, $g_0$ covers $f_0$ and $f_1$, and $g_1$ covers $f_1$. 
Next we can form the distributive lattice $H = H(I(J(L))$, which is 
isomorphic to $L$ by Theorem~\ref{isom}. 
Finally we apply the mapping $F: H \rightarrow \M_w $ defined as follows. 
First define 
$$
\X = \bigset{h \geq_T f_0 : h \T g_0 } \cup
\bigset{h \geq_T f_1 : h \T g_0 \wedge h \T g_1 }. 
$$
This has the effect that modulo $\X$ we have 
$\X \meet \{f_0\}' \equiv_w \X \meet \{g_0\}$  and 
$\X \meet \{f_1\}' \equiv_w \X \meet \{g_0, g_1\}$. 
For every $A\in H$ define 
$$
\widehat{A} = \bigset{f \in I(J(L)): f \mbox{ maximal in } A}.
$$
Finally define 
$$
F(A) = \X \meet \Meet\big\{ \{f\}' : f \in \widehat{A} \big\} \meet 
\big\{f \in I(J(L)): f \T \widehat{A} \big\}. 
$$
Here $f \T \widehat{A}$ denotes that $f \T g$ for every 
$g\in \widehat{A}$. 
We thus obtain the lattice $F(H)$ on the bottom left of the picture, with 
{\setlength\arraycolsep{2pt}
\begin{eqnarray*} 
\A &=& \X \meet \{f_0,f_1\}  \\
\B &=& \X \meet \{f_0\}' \meet \{f_1\}   \\
\C &=& \X \meet \{f_0\} \meet \{f_1\}' = 
\X \meet \{f_0, g_1\} \\
\D &=& \X \meet \{f_0\}' \meet \{f_1\}' =   
\X \meet \{f_0\}' \meet \{g_1\} = 
\X \meet \{g_0, g_1\} \\   
\E &=& \X\meet \{f_0\} \meet \{g_1\}'  \\
\F &=& \X \meet \{f_0\}' \meet \{g_1\}' = 
\X \meet \{g_0\} \meet \{g_1\}' \\
\G &=& \X \meet \{g_0\}' \meet \{g_1\}  \\
\mathcal{H} &=& \X \meet \{g_0\}' \meet \{g_1\}' 
\end{eqnarray*}
}%
Using Examples~\ref{exdiam} and~\ref{exdouble} 
one can check that $F$ is an isomorphism 
between $H$ and $F(H)$, so that the interval 
$\big[\A,\mathcal{H}\big] = \big[ F(\emptyset), F(I(J(L))) \big]$ 
is indeed isomorphic to~$L$. 
\end{example}

\noindent
We are now ready to prove:

\begin{theorem} \label{sufficient} 
Suppose that $L$ is a finite distributive lattice such that $J(L)$ is a 
an initial segment of a finite upper semilattice. 
Then $L$ is isomorphic to an interval of $\M_w$. 
\end{theorem}
\begin{proof}
We follow the procedure depicted in Figure~\ref{procedure}. 
Let $L$ be as in the hypothesis of the theorem. 
Since every finite upper semilattice with a least element is isomorphic to 
an initial segment of the Turing degrees~$\D_T$ 
(cf.\ Lerman~\cite[p156]{Lerman})
we have a finite poset $I(J(L))$ in $\D_T$ that is isomorphic to~$J(L)$. 
This means that if that if $g$ covers $f$ in $J(L)$ then 
the image of $g$ is a minimal cover of $f$ in~$\D_T$. 
Furthermore, the minimal elements of $I(J(L))$ can be chosen to be 
of minimal T-degree (so that in particular they are all noncomputable). 
Next we form the distributive lattice $H = H(I(J(L))$, which is
isomorphic to $L$ by Theorem~\ref{isom}.
Finally we define the mapping $F: H \rightarrow \M_w $ as follows.
For a given $f \in I(J(L))$ let $g_0,\ldots, g_m$ be all elements 
of $I(J(L))$ covering~$f$. 
Define
{%
$$
\renewcommand\arraystretch{2.0}
\setlength\arraycolsep{2pt}
\begin{array}{rcl} 
\X_f  &=& \bigset{h \in \omega^\omega: h >_T f \wedge h \T  g_0 
\wedge \ldots \wedge h \T  g_m },  \\
\X &=& \displaystyle 
\bigcup_{f\in I(J(L))} \X_f. 
\end{array}
$$
}%
Notice that if $f$ is maximal in $I(J(L))$ then there are no elements 
of $I(J(L))$ covering~$f$, hence $\X_f \equiv_w \{f\}'$. 
Next, for every $A\in H$ define
$$
\widehat{A} = \bigset{f \in I(J(L)): f \mbox{ maximal in } A},
$$
$$
F(A) = \X \meet \Meet\big\{ \{f\}' : f \in \widehat{A} \big\} \meet
\big\{f \in I(J(L)) : f \T \widehat{A} \big\}.
$$
Here $f \T \widehat{A}$ denotes that $f \T g$ for every
$g\in \widehat{A}$. 
By definition, $f \T \emptyset$ holds for every~$f$, so we have that 
$$
F(\emptyset) = \X \meet I(J(L)) 
\equiv_w \X \meet \big\{f \in I(J(L)) : f \mbox{ minimal in } I(J(L))\big\}.
$$
We thus obtain the lattice $F(H)$. 
Note that $H$ has $\emptyset$ as least element and 
$I(J(L))$ as largest element. 
We prove that $F$ is an isomorphism from $H$ to the interval 
$\big[F(\emptyset), F(I(J(L)))\big]\subseteq \M_w$. 
Since $H$ is isomorphic to $L$ this suffices to prove the theorem. 

{\em $F$ is injective.} 
Suppose that $A, B\in H$. Note that since $A$ and $B$ are downwards closed, 
$A = B $ if and only if $\wh{A} = \wh{B}$. 
So it suffices to show that if $\wh{A} \not\subseteq \wh{B}$ then 
$F(A) \not\equiv_w F(B)$. Suppose that $f \in \wh{A} - \wh{B}$. 
Since $f \in I(J(L))$ we have $\{f\} \not\geq_w \X$. 
Since $f\in \wh{A}$ we also have $\{f\} \not\geq_w F(A)$ because 
$\wh{A}$ is an antichain. 
If $f \T \wh{B}$ then $\{f\}\geq_w F(B)$ so in this case $F(A) \not\equiv_w F(B)$. 
If $f \nT \wh{B}$ then there is $g\in\wh{B}$ with either $f\geq_T g$ or $g>_T f$. 
In the first case we have $f>_T g$ (because $g$ is in $\wh{B}$ and $f$ is not), 
hence $\{f\} \geq_w F(B)$, and again we can conclude that $F(A) \not\equiv_w F(B)$.
In the second case, since $g>_T f \in \wh{A}$ we have $\{g\} \geq_w F(A)$, but 
$\{g\} \not\geq_w F(B)$ because $g\in\wh{B}$ and $\wh{B}$ is an antichain, 
so again $F(A) \not\equiv_w F(B)$.

{\em $F$ is monotone.} 
We claim that $A\subseteq B$ implies that $F(A) \leq_w F(B)$. 
Suppose that $A\subseteq B$ and that $h \in F(B)$. 
We prove that $\{h\} \geq_w F(A)$. We have the following three 
cases, corresponding to the three components of $F(B)$: 
\begin{itemize}

\item If $h \in \X$ then we are immediately done. 

\item If $h \in I(J(L))$, $h\T \wh{B}$ then we have one of 
the following three options:  
\begin{itemize}
\item $h \T \wh{A}$. In this case we are done immediately.  

\item $\exists g \in \wh{A} \;\; h\geq_T g$. In this case we cannot have 
$h \in \wh{A}$ because $A\subseteq B$ and $h\T \wh{B}$, so we 
have $h>_T g$ and hence $\{h\} \geq_w F(A)$. 

\item $\exists g \in \wh{A} \;\; h\leq_T g$. This case cannot 
occur because $A$ is downwards closed, hence $h$ would be in $A$, 
hence in $B$, contradicting $h\T \wh{B}$. 
\end{itemize}

\item $h>_T f$ for some $f \in \wh{B}$. 
When $f \in \wh{A}$ we are done. If $f \notin \wh{A}$ then 
since $A\subseteq B$ and $f$ is maximal in $B$ we have $f\notin A$, 
so either $f \T \wh{A}$, in which case we are done or 
$\exists g \in \wh{A} \; f\geq_T g$, in which case $f >_T g$ 
since $f \notin \wh{A}$, and again we are done. 

\end{itemize}

$F(A\cap B) \equiv_w F(A) \meet F(B)$:   
By monotonicity of $F$ we have $F(A\cap B) \leq_w F(A), F(B)$, hence also 
$F(A\cap B) \leq_w F(A)\meet F(B)$.  
For the other direction $\geq_w$, suppose that $h \in F(A\cap B)$. 
We consider the three cases corresponding to the three components of 
$F(A\cap B)$. 
\begin{itemize}

\item If $h \in \X$ then we immediately have that $\{h\} \geq_w F(A), F(B)$. 

\item Suppose that $h \in F(A\cap B)$ because 
$h >_T f $ for some $f\in \wh{A\cap B}$. 
If $f \in \wh{A}$ or $f \in \wh{B}$ then 
$\{f\}'\geq_w F(A)$ or $\{f\}'\geq_w F(B)$, hence we are done. 
If $\{h\} \geq_w \{ g \in I(J(L)) : g \T \wh{A} \}$ or 
$\{h\} \geq_w \{ g \in I(J(L)) : g \T \wh{B} \}$ then we are also done. 
Otherwise, in particular both $h \nT \wh{A}$ and $h \nT \wh{B}$, say that 
$g_0 \in \wh{A}$ and $g_1 \in \wh{B}$ are such that 
$h \nT g_0$ and $h \nT g_1$. 
It is impossible that $h \leq_T g_0, g_1$ because then 
(because $A$, $B$ downwards closed) 
$h \in A\cap B$, contradicting $f\in \wh{A\cap B}$. 
So at least one of $g_0 <_T h$ and $g_1 <_T h$ must hold. 
But in the first case we have $\{h\} \geq_w F(A)$ 
and in the second $\{h\} \geq_w F(B)$. 

\item Finally suppose that $h \T \wh{A\cap B}$.
When $h \T \wh{A}$ or $h \T \wh{B}$ then we are done, 
so suppose that neither of these hold, say 
$h \nT f \in \wh{A}$ and $h \nT g \in \wh{B}$. 
When either $f$ or $g$ is in $\wh{A\cap B}$ then 
$h \nT \wh{A \cap B}$ contrary to assumption, so we have 
that $f$, $g \notin \wh{A\cap B}$. 
When $f \geq_T g$ then $g \in A\cap B$, and because $g \notin \wh{A\cap B}$ 
there is then $h \in \wh{A\cap B}$ with $h >_T g$, contradicting 
$g \in \wh{B}$. Likewise, $ g\geq_T f$ is impossible, so we have $f \T g$. 
Hence either $h >_T f,g$ or $h <_T f,g$. 
In the latter case $h \in A \cap B$, contradicting $h \T \wh{A\cap B}$, 
and in the former case we have $\{h\} \geq_w F(A), F(B)$. 

\end{itemize}
Hence every 
$h \in F(A\cap B)$ computes an element of either $F(A)$ or $F(B)$.

$F(A\cup B) \equiv_w F(A) \join F(B)$:  
By monotonicity of $F$ we have $F(A), F(B) \leq_w F(A\cup B)$, hence also 
$F(A)\join F(B) \leq_w F(A\cup B)$. 
For the other direction, suppose that $\{h\}\geq_w F(A), F(B)$. 
We prove that $\{h\} \geq_w F(A\cup B)$. If $\{h\} \geq_w \X$ we are immediately 
done, so assume that $\{h\} \not\geq_w \X$. 
We have to prove that either 
\begin{eqnarray}
\exists f \in \wh{A\cup B} \;\; h >_T f  \label{done1} \\ 
\makebox[0mm][r]{\makebox[10mm][l]{ or }} \exists f \T \wh{A\cup B} \;\; h \geq_T f.  \label{done2} 
\end{eqnarray}
We have the following cases, corresponding to the four remaining ways 
in which $h$ can be above the components of $F(A)$ and $F(B)$ that are 
different from~$\X$: 
\begin{itemize}   
 
\item 
$\{h\} \geq_w \{f \in I(J(L)): f \T \wh{A}\}, \{f \in I(J(L)): f \T \wh{B}\}$. 
Suppose that $h \geq_T f_0 \oplus f_1$, $f_0 \T \wh{A}$ and $f_1 \T \wh{B}$. 
Note that it is not possible that $f_0 \oplus f_1 \leq_T g$ for some 
$g \in \wh{A\cup B}$ because then $f_0 \oplus f_1$ would be 
below an element of either $\wh{A}$ or $\wh{B}$, contradicting that 
$f_0 \T \wh{A}$ and $f_1 \T \wh{B}$. 
So we either have $f_0 \oplus f_1 >_T g$  for some $g \in \wh{A\cup B}$ 
or $f_0 \oplus f_1  \T \wh{A\cup B}$. 
In the first case we have $h >_T g$ and we are done by way of~(\ref{done1}).
In the second case we are done by way of~(\ref{done2}).

\item $\{h\} \geq_w \{f \in I(J(L)): f \T \wh{A}\}$ and 
$h >_T f_0$ for some $f_0 \in \wh{B}$. 
If there is such $f_0$ with $f_0 \in \wh{A\cup B}$ 
then we are done by way of~(\ref{done1}), 
so assume without loss of generality that $f_0 \in \wh{B} - \wh{A\cup B}$. 
If $h \T \wh{A\cup B}$ then we are done by~(\ref{done2}), so 
assume that $h \nT \wh{A\cup B}$, say $f \in \wh{A\cup B}$, $h\nT f$. 
We have one of the following two cases: 
\begin{itemize}

\item $h\leq_T f$. In this case $f \in \wh{B}$ is impossible by $f_0 \in \wh{B}$, 
so we must have $f\in \wh{A}$. But this contradicts
$\{h\} \geq_w \{f \in I(J(L)): f \T \wh{A}\}$. 

\item $h>_T f$. In this case we are done by way of (\ref{done1}). 

\end{itemize}

\item $\{h\} \geq_w \{f \in I(J(L)): f \T \wh{B}\}$ and $h >_T f_0$ for 
some $f_0 \in \wh{A}$. This is completely symmetric to the previous case.  

\item $h >_T f_0, f_1$ for some $f_0 \in \wh{A}$ and $f_1 \in \wh{B}$. 
If either $f_0$ or $f_1$ is in $\wh{A\cup B}$ then we are done by way of~(\ref{done1}). 
Otherwise $f_0, f_1 \notin\wh{A\cup B}$. 
When $h \T \wh{A\cup B}$ we are done by~(\ref{done2}), so 
assume that $h \nT \wh{A\cup B}$, say $h\nT f$ with $f \in \wh{A\cup B}$. 
Because $f_0, f_1 \leq_T h$ it is impossible that $h \leq_T f$, for then 
either $f_0$ would not be in $\wh{A}$ or $f_1$ would not be in $\wh{B}$. 
Hence $h >_T f$, and we are done by way of~(\ref{done1}).  

\end{itemize}

{\em $F$ is surjective.}
It remains to check that if $\C \in \big[F(\emptyset), F(I(J(L)))\big]$ then 
there is $A\in H$ such that $F(A) = \C$. 
To this end, let $A$ be a maximal subset of $I(J(L))$ such that 
$\C \geq_w F(A)$. We claim that then also $\C \leq_w F(A)$. 
Namely we have $\C\leq_w \X$ because $\C\leq_w F(I(J(L))$. 
As for the other components of $F(A)$, 
suppose that $f \in I(J(L))$, $f \T \widehat{A}$ and suppose that 
$\C$ contains no element of degree $\deg_T(f)$. 
Then $\C \geq_w F(A\cup \{f\})$, for the elements of $\C$ that are 
mapped to $f$ in the reduction $\C \geq_w F(A)$ are all $>_T f$. 
But $\C \geq_w F(A\cup \{f\})$ contradicts the maximality of $A$. 
It follows that $\C\leq_w \big\{f \in I(J(L)) : f \T \widehat{A} \big\}$.
We also have $\C\leq_w \Meet\big\{ \{f\}' : f \in \widehat{A} \big\}$. 
Namely suppose not, that is, suppose there is $f\in\widehat{A}$ 
such that $\C\not\leq_w \{f\}'$. 
Such $f$ cannot be maximal in $I(J(L))$ because $\C\leq_w F(I(J(L))$.  
Hence the set $\{g_0,\ldots, g_m\}$ of  all elements of $I(J(L))$ covering $f$ 
is nonempty. We have $\X \meet \{f\}' \equiv_w \X \meet \{g_0,\ldots, g_m\}$ 
because $\X_f \subseteq \X$ and because 
of the minimality of the $g_i$ over~$f$. 
But then $\C \geq_w F(A\cup \{g_0,\ldots, g_m\})$, contradicting the 
maximality of $A$. 
We have thus proved that $\C\equiv_w F(A)$. 
This concludes the proof of the surjectivity of $F$ and of the theorem. 
\end{proof}

\begin{theorem} \label{converse}
If $J(L)$ is not an initial segment of a finite upper semilattice 
then $L$ is not isomorphic to an interval in $\M_w$. 
\end{theorem}
\begin{proof}
We start by arguing as in Theorem~\ref{noccc}.
If $J(L)$ is not an initial segment of a finite upper semilattice then
then it must contain two incomparable elements $x_0$ and $x_1$ without 
a least upper bound. If $x_0$ and $x_1$ would have no or only one upper 
bound in $J(L)$ then $J(L)$ could be consistently extended to 
an upper semilattice, so there must be at least {\em two\/} 
incomparable minimal upper bounds $y_0$ and $y_1$ for both $x_0$ and $x_1$. 
But then the poset $\{x_0,x_1,y_0,y_1\}$ is isomorphic to the 
configuration in Figure~\ref{JL}. 
(Note that the intervals between the $x$'s and the $y$'s need not be 
empty though.)    
Consider the subinterval $[X,Y] \subseteq H(J(L))$, where 
{\setlength\arraycolsep{2pt}
\begin{eqnarray*} 
X &=& \bigset{x\in J(L): x< x_0 \vee x< x_1},  \\
Y &=& \bigset{x\in J(L): x\leq y_0 \vee x\leq y_1}.
\end{eqnarray*}
}%
Suppose that $[\X,\Y]\subseteq \M_w$ is isomorphic to $[X,Y]$
From this assumption we will derive a contradiction. 

The lattice $[X,Y]$ starts and ends with a diamond,  namely 
$[X,Y]$ has at the top the diamond with top $Y$ and bottom $Y -\{y_0,y_1\}$ 
and at the bottom the diamond with top $X\cup \{x_0,x_1\}$ and bottom $X$.  
Hence, because by Example~\ref{exdiam} there is only one way to implement the 
diamond in $\M_w$, we can argue as before that there are 
$f_0, f_1, g_0, g_1 \in \omega^\omega$ with $f_0, f_1 \not\geq_w \X$, $f_0\T f_1$, 
and $g_0,g_1 \not\geq_w \Y$, $g_0\T g_1$, such that 
{\setlength\arraycolsep{2pt}
\begin{eqnarray*} 
\X &\equiv_w& \Y \meet \{f_0,f_1\},  \\
\Y \meet \{f_0\}' \meet \{f_1\}'&\leq_w& \Y \meet \{g_0,g_1\},  \\
\Y &\equiv_w& \Y \meet \{g_0\}' \meet \{g_1\}'.
\end{eqnarray*}
}%
We now have two cases: 

{\em Case 1: $g_0,g_1 >_T f_0\oplus f_1$.} 
Note that the interval 
\begin{equation} \label{interval1}
\big[X\cup\{x_0,x_1\}, Y -\{y_0,y_1\}\big]\subseteq [X,Y].
\end{equation}
under the assumed isomorphism of $[X,Y]$ with $[\X,\Y]$ must 
correspond to the interval 
\begin{equation} \label{interval2}
\big[ \Y \meet \{f_0\}' \meet \{f_1\}', \Y \meet \{g_0,g_1\} \big] \subseteq [\X,\Y].
\end{equation}
Because $f_0\oplus f_1 \notin \Y$ (since $g_0,g_1 \not\geq_w \Y$)
we have by Dyments Theorem~\ref{DymentMw} that the subinterval 
$\big(\Y \meet \{f_0\oplus f_1\}, \Y\meet \{f_0\oplus f_1\}'\big)$ is empty. 
Hence we see that the element 
$\Y \meet \{f_0\oplus f_1\}'$ in the interval~(\ref{interval2}) 
is join-irreducible in $\M_w$, and hence also in $[\X,\Y]$. 
We derive a contradiction by showing that the interval~(\ref{interval1}) 
contains only elements that are join-reducible in $[X,Y]$. 
The elements of this interval are all supersets $A \supseteq X\cup\{x_0,x_1\}$ 
that because of the special form of $\{x_0,x_1,y_0,y_1\}$ contain 
no elements both above $x_0$ and $x_1$. 
Hence every such $A$ can always be split into a nonempty part $A_0$ of elements 
above $x_0$ and a nonempty part $A_1$ of elements above $x_1$. 
In particular every element of~(\ref{interval1}) is join-reducible in $[X,Y]$. 

{\em Case 2.} If Case 1.\ does not obtain there must be at least one 
$f$ incomparable to a $g$, say that $f_0 \T g_1$. 
It is clear that $[X,Y]$ is double diamond-like.
We derive a contradiction by showing that 
$[\X,\Y]$ is not double diamond-like. 

Let $\HH$ consist of $f_0$ plus 
the set of all elements $h$ with $f_1 \leq_T h \leq_T g_0$ that 
are minimal with respect to the property $h\not\leq_T g_1$. 
(Note that there can be only finitely many $h$ in $[f_1,g_0]$ since otherwise 
$[\X,\Y]$ would be infinite, in which case we would have reached a contradiction 
right away.) 
Then we have 
$$
\X \equiv_w \Y \meet \{f_0,f_1\} \equiv_w 
\big(\Y \meet \HH \meet \{g_1\}'\big) \meet
\big(\Y \meet \{f_0\}' \meet \{f_1\}\big).
$$
We prove that the two elements on the right hand side are join-irreducible 
in $[\X,\Y]$, 
and that one is maximal with this property. Since $\X$ is the~0 of the 
lattice $[\X,\Y]$ this proves that this interval is not double diamond-like. 

{\em $\Y \meet \{f_0\}' \meet \{f_1\}$ is join-irreducible:} 
When $\A \in [\X,\Y]$ and $\A <_w \Y\meet \{f_0\}' \meet \{f_1\}$
then $\A \leq_w \Y\meet \{f_0\} \meet \{f_1\}$: 
We have $\X \equiv_w \Y \meet \{f_0\} \meet \{f_1\} \leq_w \A$, 
and if in this reduction all elements of $\A$ that are mapped to $f_0$ 
are strictly above it then $\Y \meet \{f_0\}'\meet \{f_1\} \leq_w \A$, 
contradicting the assumption. 
Hence $\A$ contains an element of $\deg_T(f_0)$, and because 
$\A \leq_w \Y \meet \{f_1\}$ we then have 
$\A \leq_w \Y \meet \{f_0\} \meet \{f_1\}$.

{\em $\Y \meet \HH \meet \{g_1\}'$ is a 
maximal join-irreducible element of $[\X,\Y]$:} 
Join-reducibility is seen with an argument similar to the previous one: 
When $\A \in [\X,\Y]$ and 
$\A <_w \Y \meet \HH \meet \{g_1\}'$ then 
$\A \leq_w \Y \meet \HH \meet \{g_1\}$. Namely, 
we clearly have $\A \leq_w \Y \meet \HH$.
If $\A$ would not have an element of degree $\deg_T(g_1)$ then 
contrary to assumption we would have $\A \geq_w \Y \meet \HH \meet \{g_1\}'$:
By definition of $\HH$ and because $\A \geq_w \X \equiv_w \Y\meet\{f_0,f_1\}$, 
any element in $\A$ that is not below $g_1$ can be mapped 
to either $\Y$ or $\HH$. Any element in $\A$ that can be mapped to $g_1$ is 
actually strictly above $g_1$, so can be mapped to $\{g_1\}'$ by the 
identity. Hence we also have $\A \leq_w \{g_1\}$. 

To see the maximality: Suppose $\A \in [\X,\Y]$ is such that 
$\A >_w \Y \meet \HH \meet \{g_1\}'$. 
Then $\A$ is of the form $\A = \Y \meet \mathcal{K} \meet \{g_1\}'$, 
with $\{g_0\}' \geq_w \mathcal{K} >_w \HH$. But then it is easy to check that 
$$
\A \equiv_w \big(\Y \meet \HH \meet \{g_1\}'\big) \join 
\big(\Y \meet \mathcal{K} \meet \{g_1\} \big), 
$$
using that $\HH \oplus \{g_1\} \subseteq \{g_1\}'$. 
The two components on the right hand side are w-incomparable and in $[\X,\Y]$, 
so $\A$ is join-reducible in $[\X,\Y]$. This concludes Case 2.

Summarizing, 
we have seen that $H(J(L))$ contains a subinterval $[X,Y]$ that 
cannot be an interval in $\M_w$. 
Thus $H(J(L))$, and hence by  by Theorem~\ref{isom} also $L$, 
cannot be an interval of $\M_w$.
\end{proof}

\noindent
By combining the above results we obtain the following 
characterization of the finite intervals of $\M_w$:

\begin{theorem} \label{main}
For any finite distributive lattice $L$ the following are equivalent:
\begin{enumerate}

\item $L$ is isomorphic to an interval in $\M_w$, 

\item $J(L)$ is an initial segment of a finite upper semilattice, 

\item $L$ has no double diamond-like lattice as a subinterval. 

\end{enumerate} 
\end{theorem}
\begin{proof}
Item~2.\ and 3.\ are equivalent by Theorem~\ref{noccc}. 
They imply item~1.\ by  Theorem~\ref{sufficient}. 
Conversely, 1.\ implies~2.\ by Theorem~\ref{converse}. 
\end{proof}

\begin{corollary} \label{initials}
A finite distributive lattice is an initial segment of $\M_w$ 
if and only if it has no double diamond-like subinterval 
and it has a meet-irreducible~0.  
\end{corollary}
\begin{proof}
We can extend the definition of the mapping $F$ in the 
proof of Theorem~\ref{main} as follows. 
Define $\X_f$ as before and let 
{%
$$
\renewcommand\arraystretch{2.0}
\setlength\arraycolsep{2pt}
\begin{array}{rcl} 
\X_0  &=& \bigset{h \in \omega^\omega : 
h \T  f \mbox{ for all $f$ minimal in } I(J(L)) },  \\
\X &=& \displaystyle \X_0 \cup \bigcup_{f\in I(J(L))} \X_f. 
\end{array}
$$
}%
Then for every $A \in H$ define $F(A)$ as before, using this 
new definition of $\X$. 
This addition does not change anything in the proof of Theorem~\ref{sufficient}, 
but now we have that $F(\emptyset) \equiv_w 0'$, as is easily checked, 
using that we chose the minimal elements of $I(J(L))$ of minimal T-degree. 
Thus we obtain that a finite distributive lattice has no double diamond-like 
subinterval if and only if it is isomorphic to an interval of the form 
$[\0',\bfA]$ in~$\M_w$. From this the corollary follows immediately. 
\end{proof}

\begin{acknowledgments}
The author thanks Andrea Sorbi for helpful discussions. 
\end{acknowledgments}

\end{document}